\documentclass{amsart}

\usepackage[english]{babel}
\usepackage{amssymb}
\usepackage{mathptmx}
\usepackage[numbers]{natbib}
\usepackage{graphicx}
\usepackage{epsf, subfigure, verbatim}

\theoremstyle{remark}

\DeclareMathOperator*{\asinh}{arcsinh}

\begin{document}
\title[Estimation of first passage time densities]{{\itshape Estimation of first passage time densities of diffusions processes\break through time-varying boundaries} \break}
\author[Imene Allab and Francois Watier]{}
\maketitle
\begin{center}
{\large Imene Allab and Francois Watier}\\
\smallskip
\smallskip
{\small {\em Department of Mathematics, Universit\'e du
Qu\'ebec \`a Montr\'eal, Montr\'eal, Canada}}
\end{center}
\bigskip

\begin{abstract}
In this paper, we develop a Monte Carlo based algorithm for estimating the FPT density of a time-homogeneous SDE through a time-dependent frontier. We consider Brownian bridges as well as localized Daniels curve approximations to obtain tractable
estimations of the FPT probability between successive points of a simulated path of the process.
Under mild assumptions, a (unique) Daniels curve local approximation can easily be obtained by explicitly solving a non-linear system of equations.\medskip



\end{abstract}

\section{Introduction}

Let $X$ be a time homogeneous diffusion process which is the unique
(strong) solution of the following stochastic differential equation
:
\begin{equation}\label{diffX}
dX\left(t\right)=\mu\left(X\left(t\right)\right)dt+\sigma\left(X\left(t\right)\right)dW\left(t\right),\:X\left(0\right)=x_{0}
\end{equation}

if $S$ is a time dependent boundary, we are interested in estimating
either the pdf or cdf of the first passage time (FPT) of the diffusion process
through this boundary that is we will study the following random variable :
\[
\tau_{S}=\inf\left\{ t>0|X\left(t\right)=S\left(t\right)\right\}
\]

In general, there is no explicit expression for the first passage-time
density of a diffusion process through a time-varying boundary. To
this date, only a few specific cases provide closed formed formulas for example
when the process is gaussian and the boundary is of a Daniels' curve
type. Thus, we mainly rely on simulation techniques
to estimate this density in a general setting.
\medskip{}

The main goal of this work is to develop a computationally efficient algorithm that will provide reliable FPT density estimates. The paper is organized as follows. In section 2, we review existing techniques followed by the mathematical foundations leading to a novel algorithm. Finally, section 3 is devoted to various examples enabling us to evaluate the algorithm's performance.

\section{Monte Carlo simulation estimation}

\medskip{}

This is the simplest and best-known approach based on the law of large numbers. After fixing
a time interval, basically we divide the latter into smaller ones,
simulate a path of the process along those time points and, if it
occurs, note the subinterval where the first upcrossing occurs. Generally,
the midpoint of this subinterval forms the estimated first passage
time of this simulated path. We repeat the process a large number
of time to construct a pdf or cdf estimate of this stopping time.
\medskip{}

Consider a Brownian motion $W$, a linear boundary $S\left(t\right)=\alpha+\beta t$, then for $\alpha>0$ and $T>0$, from standard theory the first passage time probability has an explicit form given by
\[
P\left(\tau_{S}\leqslant T\right)=\Phi\left(-\frac{\alpha+\beta T}{\sqrt{T}}\right)+e^{-2\alpha\beta}\Phi\left(\frac{-\alpha+\beta T}{\sqrt{T}}\right)
\]
where $\Phi$ denotes the cdf of a standard normal distribution.\\

Setting $\alpha=0.5$, $\beta=0.2$ and $T=1$, table 1 gives us estimates of the FPT probability with various number
of simulated paths $N$ and time-step discretization $\Delta t$. Clearly, even in
a simple case as this one, in order to have a suitable estimation of the true value $P\left(\tau_{S}\leqslant T\right)=0.5548$
we have to rely on a large number of paths and a very fine partition of the time interval.
\medskip{}

\begin{figure}[h]
  \centering
  \includegraphics[scale=0.6]{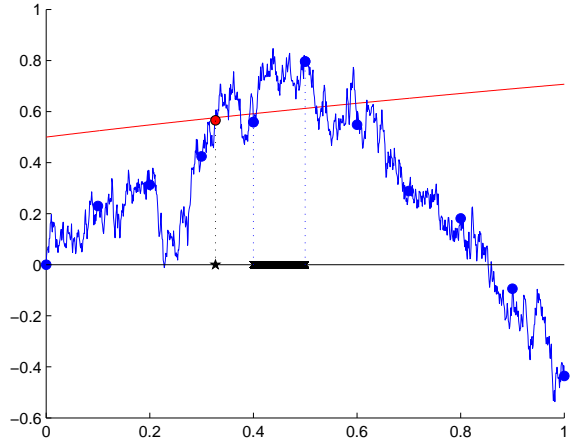}\\
  \caption{Undetected prior upcrossing through a basic Monte Carlo path simulation}\label{onepath}
\end{figure}

\begin{table}[h]
\begin{center}
\caption{Monte Carlo estimates of the FPT probability of a Brownian motion through the boundary $S\left(t\right)=0.5+0.2 t$}
\begin{tabular}{|c|c|c|c|}
\hline
$P\left(\tau_{S}\leqslant 1.0\right)$  & $\Delta t=10^{-2}$  & $\Delta t=10^{-3}$  & $\Delta t=10^{-4}$ \\
\hline
\hline
$N=10^{4}$  & 0.5096  & 0.5407  & 0.5487\\
\hline
$N=10^{5}$  & 0.5102  & 0.5390  & 0.5506\\
\hline
$N=10^{6}$  & 0.5122  & 0.5412  & 0.5504\\
\hline
\end{tabular}
\par\end{center}
\end{table}
Another drawback of the crude Monte Carlo approach is that it
tends to overestimates the true value of the first passage-time since
an upcrossing may occur earlier in between simulated points of a complete
path as illustrated in figure \ref{onepath}.

\subsection{Monte Carlo approach with intertemporal upcrossing simulations}

Instead of continuously repeating the whole Monte Carlo procedure with an even
finer interval partition to obtain better estimates, let us see how
one could improve on the initial estimates without discarding the
simulated paths.\medskip{}

An astute idea that have been put forward by several authors, is to ideally obtain the probability law of an upcrossing
between simulated points, thus if $p_{k}$ is the probability of an
upcrossing in the time interval $\left[t_{k},t_{k+1}\right]$, then
one would simply generate a value $U_{k}$ taken from a uniform random
variable on $\left[0,1\right]$ and assert that there is an upcrossing
if $U_{k}\leqslant p_{k}$. Since the exact FPT probability of a diffusion
bridge will more than often not be available, we need to consider
an adequate estimation of this probability.

\subsubsection{Diffusion bridge approximation}

\medskip{}

For each subinterval, one could consider simulating paths of approximate
tied-down processes as proposed by Giraudo, Sacerdote and Zucca \cite{GSZ2001} where they basically used a Kloden
Platen approximation scheme with order of convergence 1.5, or we could
make use of more recent results from Lin, Chen and Mykland \cite{LCM2010} or S\o{}rensen and Bladt \cite{SB2013} to improve on the
FPT probability estimates. Although all of these may constitute adequate approximations
of the true FPT probability it may prove costly in computation time
since these methods tantamounts to generating numerous simulations of bridge paths on successive subintervals for each of the original sample paths.\\

Another alternative, as first proposed in Strittmatter \cite{S1987}, is to consider that
for a small enough interval, the diffusion part of the process should
remain fairly constant and then consider a Brownian bridge approximation
of the diffusion bridge and exploit known results on the FPT of Brownian
bridges. For example Giraudo and Sacerdote \cite{GS1999} considered solving numerically the Volterra type integral equation linked to the
generalized Brownian bridge FPT probability through a general time varying
boundaries. Since a certain number of iterations may be needed to obtain adequate solutions of integral equations specific to each sample paths and successive subinterval this may sensibly increase computing time.

\medskip{}

Finally, it is worth mentioning that all of the above methods could, in many cases, be
improved significantly as far as accuracy is considered by first applying
a Lamperti transform on both the original process and the frontier
as described for example in Iacus \cite{I2008}.

\medskip{}

Indeed , define
\begin{equation}\label{Lamp}
F\left(x\right)=\int^{x}\frac{1}{\sigma\left(u\right)}du
\end{equation}
and apply it on the original process and the time-varying boundary.
Assuming that $F$ is one-to-one, then the original problem is equivalent
to finding the FPT density of
\[
\tau_{S^{*}}=\inf\left\{ t>0|Y\left(t\right)=S^{*}\left(t\right)\right\}
\]

where the new boundary is given by $S^{*}\left(t\right)=F\left(S\left(t\right)\right)$
and, by It\^o's formula, the diffusion process $Y$ follows the dynamic
\begin{equation}\label{diffY}
dY\left(t\right)=b\left(Y\left(t\right)\right)dt+dW\left(t\right),\;Y\left(0\right)=F\left(x_0\right)
\end{equation}

where
\[
b\left(t,Y\left(t\right)\right)=\frac{\mu\left(F^{-1}\left(Y\left(t\right)\right)\right)}{\sigma\left(F^{-1}\left(Y\left(t\right)\right)\right)}-\frac{1}{2}\frac{\partial\sigma}{\partial x}\left(F^{-1}\left(Y\left(t\right)\right)\right)
\]
Since the diffusion part of the process $Y$ is constant, then the
simple Brownian bridge will constitute a good approximation of the
diffusion bridge.

\medskip{}

\subsubsection{Diffusion bridge approximation with local boundary approximation}

In our approach, while still considering a Brownian bridge approximation
of the diffusion bridges after a Lamperti transform as described previously, we propose to consider
localized Daniels curve approximation of the time-varying boundary.
Since explicit formulas of the first passage time probability are
available in this case, one would readily get an adequate approximation
of the true probability $p_{k}$. Furthermore, under mild assumptions,
a (unique) Daniels curve approximation can easily be obtained by simply
taking the endpoints of the segment and the value at midpoint (or
another point of our choosing). Indeed, we will show that it leads
to consider a non-linear system of three equations that can be explicitly solved.

\medskip{}

Before describing our algorithm, we will need the following
key results :

\medskip{}

\textbf{Proposition 1}. Consider a Brownian bridge $W^{0,\Delta t}$ define
on an time interval $\left[0,\Delta t\right]$ and $S$ a Daniels curve defined
by
\begin{equation}\label{Dan}
S\left(t\right)=\frac{\alpha}{2}-\frac{t}{\alpha}\ln\left[\frac{\beta+\sqrt{\beta^{2}+4\gamma e^{-\alpha^{2}/t}}}{2}\right]
\end{equation}
where $\alpha,\beta>0$, $\gamma\in R$ and $\lim_{t\rightarrow \Delta t}\beta^{2}+4\gamma e^{-\alpha^{2}/t}>0$,
if $\tau_{S}=\inf\left\{ 0\leqslant t\leqslant \Delta t|W^{0,\Delta t}(t)=S(t)\right\} $
then
\[
P\left(\tau_{S}\leqslant \Delta t\right)=\beta e^{-\frac{\alpha^{2}}{2\Delta t}}+\gamma e^{-\frac{2\alpha^{2}}{\Delta t}}
\]

\emph{Proof. } Apply Therom 3.4 of Di Nardo et al. \cite{DN2001}.

\medskip{}

\textbf{Proposition 2}. Let $\left[0,\Delta t\right]$ be a time interval,
consider the points $\left(0,a\right)$, $\left(\Delta t/2,b\right)$,$\left(\Delta t,c\right)$
and set $A=e^{\frac{2a^{2}}{\Delta t}}$ $B=e^{-\frac{4\alpha b}{\Delta t}}$et
$C=e^{-\frac{2ac}{\Delta t}}$. If
\begin{equation}\label{cond1}
\frac{1}{A^{2}}<\frac{B}{C}<\frac{1}{A}+\frac{\sqrt{A^{2}-1}}{A^{2}}
\end{equation}

then there is a unique Daniels curve (\ref{Dan}) passing through the three points
with parameters
\begin{equation}\label{param}
\alpha=2a,\:\beta=\frac{A\left(A^{4}B^{2}-C^{2}\right)}{A^{3}B-C},\:\gamma=A^{4}C^{2}-\beta A^{3}C
\end{equation}

\emph{Proof. }The set of points generate the following non-linear
system of equations :
\begin{eqnarray}\label{eq1}
\nonumber\frac{\alpha}{2} & = & a\\
\frac{\alpha}{2}-\frac{\Delta t}{2\alpha}\ln\left[\frac{\beta+\sqrt{\beta^{2}+4\gamma e^{-2\alpha^{2}/\Delta t}}}{2}\right] & = & b\\
\nonumber\frac{\alpha}{2}-\frac{\Delta t}{\alpha}\ln\left[\frac{\beta+\sqrt{\beta^{2}+4\gamma e^{-\alpha^{2}/\Delta t}}}{2}\right] & = & c
\end{eqnarray}

obviously the first equation gives $\alpha=2a>0$, while simple algebraic
manipulations on the last two equations lead us to solve the following
linear system
\begin{eqnarray*}
\beta e^{\frac{\alpha^{2}-2\alpha b}{\Delta t}}+\gamma e^{-\frac{2\alpha^{2}}{\Delta t}} & = & e^{\frac{2\alpha^{2}-4\alpha b}{\Delta t}}\\
\beta e^{\frac{\alpha^{2}-\alpha b}{\Delta t}}+\gamma e^{-\frac{\alpha^{2}}{\Delta t}} & = & e^{\frac{\alpha^{2}-2\alpha b}{\Delta t}}
\end{eqnarray*}

which can be rewritten in the form
\begin{eqnarray*}
\beta A^{6}B+\gamma & = & A^{8}B^{2}\\
\beta A^{3}C+\gamma & = & A^{4}C^{2}
\end{eqnarray*}

since $A^{6}B-A^{3}C=A^{3}\left(A^{3}B-C\right)>0$ then there exists
a unique solution given by :
\begin{eqnarray*}
\beta & = & \frac{A\left(A^{4}B^{2}-C^{2}\right)}{A^{3}B-C}\\
\gamma & = & A^{4}C^{2}-\beta A^{3}C\\
 & = & A^{7}BC\left[\frac{C-AB}{A^{3}B-C}\right]
\end{eqnarray*}

this would constitute the solution to the original system provided
that $\beta>0$ and $\lim_{t\rightarrow \Delta t}\beta^{2}+4\gamma e^{-\alpha^{2}/t}>0$.

\medskip{}

Notice first that $\frac{1}{A^{2}}<\frac{B}{C}\Rightarrow A^{4}B^{2}-C^{2}>0$
and therefore $\beta>0$ and if furthermore $\frac{B}{C}\leqslant\frac{1}{A}$
then $\gamma\geqslant0$ and clearly $\lim_{t\rightarrow \Delta t}\beta^{2}+4\gamma e^{-\alpha^{2}/t}>0$
is satisfied. So if we assume now that $\frac{B}{C}>\frac{1}{A}$
then $\gamma<0$ , thus we need to verify that $\beta^{2}+\frac{4\gamma}{A^{2}}>0$
which is the case since

\begin{eqnarray*}
\beta^{2}+\frac{4\gamma}{A^{2}} & = & \frac{A^{2}\left(A^{4}B^{2}-C^{2}\right)^{2}}{\left(A^{3}B-C\right)^{2}}+4A^{5}BC\left[\frac{C-AB}{A^{3}B-C}\right]\\
 & = & \frac{A^{2}}{\left(A^{3}B-C\right)^{2}}\left(\left(A^{4}B^{2}-C^{2}\right)^{2}+4A^{3}BC\left(C-AB\right)\left(A^{3}B-C\right)\right)\\
 & = & \frac{A^{2}}{\left(A^{3}B-C\right)^{2}}\left(A^{4}B^{2}+C^{2}-2A^{3}BC\right)^{2}
\end{eqnarray*}

The final step is to make sure that it solves the original system.
Substituting back in (\ref{eq1}), (where only positive square roots are involved), we see that is the case only if $A^{4}B^{2}+C^{2}-2A^{3}BC<0$,
or equivalently
\[
\left(\frac{B}{C}-\frac{1}{A}\right)^{2}<\frac{A^{2}-1}{A^{4}}
\]
which is verified through (\ref{cond1}).\hfill\qed\\

The FPT algorithm is described as follows :
\begin{enumerate}
 \item[Step 1] Apply the Lamperti transform (\ref{Lamp}) to the original diffusion process (\ref{diffX}) and frontier $S$ to obtain the new process (\ref{diffY}) and boundary $S^*=F\left(S\left(t\right)\right)$
  \item[Step 2] Select a time interval $\left[T_l,T_u\right]$\ and construct a partition $T_l=t_{0}<t_{1}<\ldots<t_{n}=T_u$
  \item[Step 3] Initialize FPT vector counter to $\tau:=\left\{0,\ldots,0\right\}$
  \item[Step 4] Initialize path counter to $k:=1$\\

\hspace{-1.6cm}WHILE $k$ is less than $N$ the number of desired paths DO the following :\\
\begin{enumerate}
  \item [Step 5] Simulate a path of the process $\left\{Y\left(t_{1}\right),\ldots,Y\left(t_{n}\right)\right\}$
  \item [Step 6] Initialize subinterval counter to $i:=1$\\

\hspace{-1.6cm}WHILE $i$ is less than $n$ the number of desired subintervals DO the following :\\
\begin{enumerate}
  \item [Step 7] IF $Y\left(t_{i}\right)\geqslant S^*\left(t_{i}\right)$ THEN
set $i^{th}$ FPT vector component to $\tau_{i}:=\tau_{i}+1$ and path counter to $k:=k+1$, GO TO Step 5
  \item [Step 8] Set $\Delta:=t_{i}-t_{i-1}$, $a:=S\left(t_{i-1}\right)-Y\left(t_{i-1}\right)$, $b:=S\left(\frac{t_{i-1}+t_{i}}{2}\right)-Y\left(\frac{t_{i-1}+t_{i}}{2}\right)$,\\ $c:=S\left(t_{i}\right)-Y\left(t_{i}\right)$, finally set $A$, $B$, $C$, $\alpha$, $\beta$ and $\gamma$ as in (\ref{param}) of proposition 2
  \item [Step 9] IF $\frac{1}{A^{2}}<\frac{B}{C}<\frac{1}{A}+\frac{\sqrt{A^{2}-1}}{A^{2}}$ THEN set $c_1:=\beta$, $c_2:=\gamma$\\ IF $\frac{1}{A^{2}}\geq\frac{B}{C}$ THEN set $c_1:=0$, $c_2:=A^{4}C^{2}$\\ IF $\frac{B}{C}\geq\frac{1}{A}+\frac{\sqrt{A^{2}-1}}{A^{2}}$ THEN set $c_1:=2AC$ $c_2:=-A^{4}C^{2}$
  \item [Step 10] Set probability upcrossing to $p:=c_1 e^{-\frac{\alpha^{2}}{2\Delta}}+c_2e^{-\frac{2\alpha^{2}}{\Delta}}$
  \item [Step 11] Generate a value $U$ taken from a uniform random variable
  \item [Step 12] IF $U\leq p$ THEN set $i^{th}$ FPT vector component to $\tau_{i}:=\tau_{i}+1$ and path counter to $k:=k+1$, GO TO Step 5, ELSE set $i:=i+1$, GO TO Step 7
\end{enumerate}
\end{enumerate}
\end{enumerate}
Note that step 9 includes extreme cases where the middle point of the frontier in a subinterval may not be reached by a Daniels curve, thus we use the closest curve possible.
\section{Examples}

We will focus our examples on diffusion processes which paths can
be simulated exactly. Therefore with known results on FPT density
and bounds, it will allow us to better visualize the approximation error
due essentially to the algorithm.

\medskip{}

Example 1. Consider the following Ornstein-Uhlenbeck process and time varying boundary
\begin{eqnarray*}
  dX\left(t\right) &=& \left(1.0-0.5X\left(t\right)\right)dt+dW\left(t\right),\:X\left(0\right) = 1.6\\
  S\left(t\right) &=& 2.0\left(1.0-\sinh\left(0.5t\right)\right)
\end{eqnarray*}

This diffusion process is a Gauss-Markov process and according
to Di Nardo et al. \cite{DN2001} the chosen boundary allows us to obtain an explicit FPT density
given by
$$f\left(S\left(t\right),t\right)=\frac{e^{0.5t}}{\sinh\left(0.5t\right)}\varphi_{0}\left(S\left(t\right)\right)$$

where $\varphi_{0}$ is the probability density function of the Ornstein-Uhlenbeck process starting at $X\left(0\right) = 0$.\\

\begin{figure}[h]
  \centering
  \includegraphics[scale=0.6]{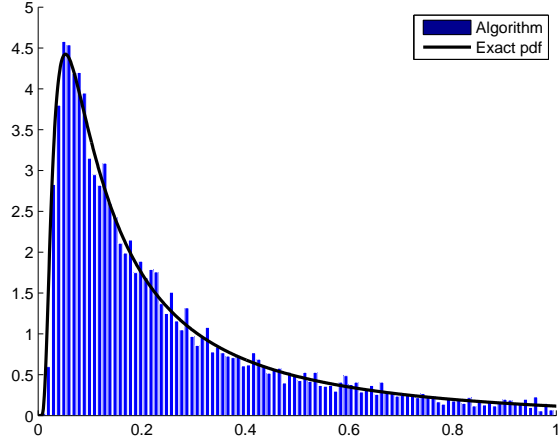}\\
  \caption{Ornstein-Uhlenbeck FPT pdf through the boundary $2.0\left(1.0-\sinh\left(0.5t\right)\right)$}\label{ex1}
\end{figure}

Figure \ref{ex1} compares the true FPT density with the empirical density
histogram obtained through our algorithm using a time step discretization
of 0.01 and 10 000 simulated paths. Furthermore, the algorithm gives us a FPT probability estimate of $0.9622$ over the whole interval compared to the true value of $0.9608$ representing a relative error of about $0.15\%$.

Example 2. Consider the following geometric Brownian process and linear boundary
\begin{eqnarray*}
  dX\left(t\right) &=& 5.0X\left(t\right)dt+2.5X\left(t\right)dW\left(t\right),\;X\left(0\right) = 0.5\\
  S\left(t\right) &=& 1.0+2.0t
\end{eqnarray*}
By applying the Lamperti transform to both the process and boundary
we obtain respectively
\begin{eqnarray*}
  dY\left(t\right) &=& 0.75dt+dW\left(t\right),\:Y\left(0\right) = 0.4\ln\left(0.5\right)\\
  S^{*}\left(t\right) &=& 0.4\ln\left(1.0+2.0t\right)
\end{eqnarray*}
As in example 1, this transformed diffusion process is also a Gauss-Markov process and,
although the new frontier does not allow an explicit FPT density,
using the deterministic algorithm in Di Nardo et al. \cite{DN2001} with a 0.01 time step discretization,
we can obtain a reliable approximation.\\

Figure \ref{ex1} compares the Di Nardo
FPT density approximation with the empirical density histogram obtained
through our algorithm using the same time step discretization with
10 000 simulated paths. In addition, the algorithm offers a FPT probability estimate of $0.8251$ over the whole interval agreeing with the actual value of $0.8258$ (a relative error of about $0.08\%$).
\begin{figure}[h]
  \centering
  \includegraphics[scale=0.6]{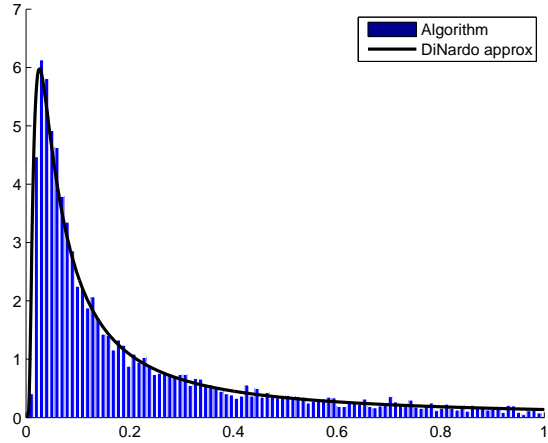}\\
  \caption{Geometric Brownian Motion FPT pdf through the boundary $1.0+2.0t$}\label{ex2}
\end{figure}

Example 3. Consider the modified Cox-Ingersoll-Ross process and linear boundary
\begin{eqnarray*}
  dX\left(t\right) &=& -0.5X\left(t\right)dt+\sqrt{1+X\left(t\right)^{2}}dW\left(t\right),\:X\left(0\right) = 0\\
  S\left(t\right) &=& 0.3+0.2t
\end{eqnarray*}
By applying the Lamperti transform to both the process and boundary
we obtain respectively
\begin{eqnarray*}
  dY\left(t\right) &=& -\tanh\left(Y\left(t\right)\right)dt+dW\left(t\right),\:Y\left(0\right) = 0\\
  S^{*}\left(t\right) &=& \asinh\left(0.3+0.2t\right)
\end{eqnarray*}
As opposed to the preceding examples, this transformed diffusion
process is not gaussian however using Beskos and Roberts' \cite{BR2005} exact algorithm we can simulate exact sample
paths. Although an explicit FPT density is not available, using results
of Downes and Borovkov \cite{DB2008} we can, in this case, obtain the following lower and upper bounds:
\begin{eqnarray*}
  f_{L}\left(S\left(t\right),t\right)&=&\frac{1}{t}\left(S\left(t\right)-\frac{0.2t}{\sqrt{1.09}}\right)\frac{e^{-0.5t}}{\cosh\left(S\left(t\right)\right)}\varphi_{W}\left(S\left(t\right)\right)\\
  f_{U}\left(S\left(t\right),t\right)&=&\frac{1}{t}\left(S\left(t\right)-\frac{0.2t}{\sqrt{1.25}}\right)\frac{e^{0.5t}}{\cosh\left(S\left(t\right)\right)}\varphi_{W}\left(S\left(t\right)\right)
\end{eqnarray*}
where $\varphi_{W}$ is the probability density function of a standard Brownian motion.\\

Figure \ref{ex1} compares the FPT bounds with the empirical density histogram obtained
through our algorithm using a 0.01 time step discretization starting
initially with 15 000 simulations and obtaining 11 768 valid paths through the exact algorithm. Moreover, the algorithm suggests a FPT probability estimate of $0.7398$ over the whole interval which lies within the values $\int_{0}^{T}f_{L}\left(S\left(t\right),t\right)dt=0.6204$ and $\int_{0}^{T}f_{U}\left(S\left(t\right),t\right)dt=0.7673$.
\begin{figure}[h]
  \centering
  \includegraphics[scale=0.6]{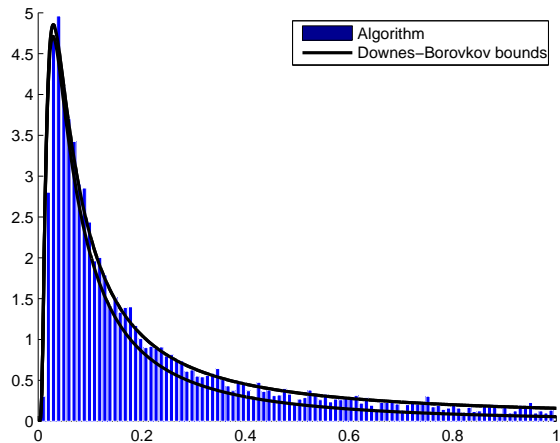}\\
  \caption{Modified CIR process FPT pdf through the boundary $0.3+0.2t$}\label{ex3}
\end{figure}
\medskip{}

\end{document}